# CROSS-VALIDATION IN NONPARAMETRIC REGRESSION WITH OUTLIERS

By Denis Heng-Yan Leung

*Singapore Management University*

A popular data-driven method for choosing the bandwidth in standard kernel regression is cross-validation. Even when there are outliers in the data, robust kernel regression can be used to estimate the unknown regression curve [*Robust and Nonlinear Time Series Analysis. Lecture Notes in Statist.* (1984) **26** 163–184]. However, under these circumstances standard cross-validation is no longer a satisfactory bandwidth selector because it is unduly influenced by extreme prediction errors caused by the existence of these outliers. A more robust method proposed here is a cross-validation method that discounts the extreme prediction errors. In large samples the robust method chooses consistent bandwidths, and the consistency of the method is practically independent of the form in which extreme prediction errors are discounted. Additionally, evaluation of the method's finite sample behavior in a simulation demonstrates that the proposed method performs favorably. This method can also be applied to other problems, for example, model selection, that require cross-validation.

**1. Introduction.** Since Nadaraya [20] and Watson [27] first proposed using the kernel method for curve estimation, there have been numerous investigations about its theory and application [9]. When there is evidence that the data may be contaminated with outliers, robust kernel regression is effective in modeling the underlying curve [4, 7].

The standard error criterion for evaluating a statistical estimator is to determine how close it is to the true parameter. In nonparametric regression, a few popular criteria are the integrated squared error, the mean integrated squared error, the average squared error and the mean average squared error. Each of these criteria gives a measure of the distance between the regression and the unknown curve "averaged" over the range of the independent variable. In practice, any one of the above criteria can be used as they are all asymptotically very similar [9]. On the other hand, the criteria are all critically dependent on the bandwidth. The bandwidth selection problem is to find the (optimal) bandwidth that is optimal in the sense that the chosen error criterion is minimized with respect to all bandwidth choices.

This paper focuses on bandwidth selection using the cross-validation method [11, 14, 21, 23]. The cross-validation method provides an estimate of the prediction error of the regression, and that, in turn, is an approximation of the error criterion. The framework of the method is straightforward. For each observation in the data, the method evaluates the prediction error using kernel regression with that observation removed from the modeling process. The optimal bandwidth is chosen to minimize the sum of squares of the prediction errors from all of the observations. Cross-validation is one of a larger class of bandwidth selection methods called data-driven methods, in which the Akaike Information Criterion [1] and Shibata's criterion [25] are examples. All these methods give asymptotically similar bandwidths [11] and, therefore, this investigation will only focus on the cross-validation method.

In standard kernel regression, where there are no outliers in the data, cross-validation has been shown to produce bandwidths that are asymptotically consistent [11, 12]. However, when there are outliers in the data, it has been demonstrated in simulations that the use of cross-validation can lead to extremely biased bandwidth estimates [17]. Arguably, the reason standard cross-validation fails, even when applied using a robust regressor [7], is because it no longer produces a reasonable estimate of the prediction error. Therefore, a robust cross-validation method may be superior. The robust cross-validation method proposed here only differs from the average squared error by a constant shift and a constant multiple; both of which are asymptotically independent of the bandwidth. Hence, the robust cross-validation rule asymptotically selects a bandwidth that optimizes the average squared error (or any other asymptotically equivalent error criterion). Similar methods have been suggested by Leung, Marriott and Wu [17] in kernel smoothing and by Ronchetti, Field and Blanchard [24] in linear model selection. Further, Wang and Scott [26] suggested an $L_1$ cross-validation in nonparametric regression. One weakness of these previous works is that they did not provide any analytical analysis of the methods. Boente, Fraiman and Meloche [2] considered a robust plug-in estimator that is based on minimizing the mean integrated squared error of the robust smoother. But that method is not fully automatic in the sense that a "pilot" bandwidth is required.

The rest of this paper is organized as follows. The method and its large sample results are given in Section 2. A simulated example and some finite sample simulation results are reported in Section 3. The method is applied to a real dataset in Section 4. A discussion of the results is given in Section 5. Proofs are given in the Appendix.

**2. Main results.** Let $(x_1, y_1), \ldots, (x_n, y_n)$ be a set of data and consider the regression of $Y$ on $X$ at the $n$ design points $x_1, \ldots, x_n$,

$$y_i = m(x_i) + \varepsilon_i, \qquad i = 1, \ldots, n, \tag{1}$$

where $m(\cdot)$ is an unknown functional of $X$ and $\{\varepsilon_i, i = 1, \ldots, n\}$ are i.i.d. random noise with distribution $F(\cdot)$. Without loss of generality, it is assumed that the range of $X$ is $[0, 1]$.

A robust smoother $\tilde{m}(\cdot)$ of $m(\cdot)$ is defined by

$$n^{-1} \sum_{i=1}^{n} \tilde{\rho}\{y_i - \tilde{m}(x_i)\} = \min,$$

where $\tilde{\rho}(\cdot)$ is an even function with bounded first derivative $\tilde{\psi}(\cdot)$. Consequently, the robust smoother can also be defined as the zero of

$$\sum_{i=1}^{n} \alpha_i(x) \tilde{\psi}(y_i - \cdot), \tag{2}$$

where $\alpha_i(x)$ is a weight function based on a kernel $K$ and the bandwidth $h$. The choice of $\tilde{\rho}$ is problem-specific. A well-known choice can be found in [16]. Some common forms of the weight function $\alpha_i(x)$ are

$$(nh)^{-1} K\left(\frac{x - x_i}{h}\right) \bigg/ \sum_{i=1}^{n} K\left(\frac{x - x_i}{h}\right) \qquad ([20, 27]),$$

$$h^{-1} \int_{t_{i-1}}^{t_i} K\left(\frac{u - x}{h}\right) du, \qquad t_0 = 0, t_n = 1, t_i = 1/2(x_i + x_{i+1})$$

$$\text{for } i = 1, \ldots, n-1 \ ([5]),$$

$$(x_i - x_{i-1}) \frac{1}{h} K\left(\frac{x - x_i}{h}\right) \qquad \text{for } i = 2, \ldots, n \ ([22]).$$

Similarly, the leave-one-out-smoother, $\tilde{m}_{-i}(\cdot)$ of $m(\cdot)$, is defined as the zero of

$$\sum_{j \neq i} \alpha_j(x) \tilde{\psi}(y_j - \cdot). \tag{3}$$

While there are a number of studies that demonstrate the effectiveness of the robust smoother, $\tilde{m}$, in estimating $m$ [4, 7], the result of directly applying the cross-validation rule, even using the robust smoother, is not as satisfactory. In fact, simulation results [17] showed that the conventional cross-validation rule,

$$CRVD(h) = n^{-1} \sum_{i=1}^{n} \{y_i - \tilde{m}_{-i}(x_i)\}^2,$$

behaves unsatisfactorily when the data is contaminated with outliers, a situation when robust smoothing is most needed. The form of the CRVD suggests the root of the problem. The CRVD is a sum of squares of the prediction errors of the smoother at each of the design points. When there are outliers, some of the prediction errors will be uncharacteristically extreme and these extreme prediction errors will affect the performance of the CRVD. Therefore, the extreme prediction errors should be discounted, as if they were extreme observations from a set of data. Based on this argument, a robust cross-validation rule can be defined as

$$(4) \qquad RCRVD(h) = n^{-1} \sum_{i=1}^{n} \rho\{y_i - \tilde{m}_{-i}(x_i)\},$$

where $\rho(\cdot)$ is a function whose role is similar to $\tilde{\rho}$. $RCRVD(h)$ can be used as a surrogate for

$$ASE(h) = n^{-1} \sum_{i=1}^{n} \{\tilde{m}(x_i) - m(x_i)\}^2$$

and

$$MASE(h) = E\{ASE(h)|x_1,\ldots,x_n\}.$$

The choice of $\rho$ and $\tilde{\rho}$ will be discussed further. It will be demonstrated that (4) is a reasonable bandwidth selector, using the following assumptions:

(A1) $F(\cdot)$ is symmetric about zero.
(A2) The kernel $K$ is symmetric about zero. Furthermore, $K$ is positive, Lipschitz continuous and satisfies

$$(5) \qquad \int K(t)\,dt = 1.$$

(A3) The function $m:[0,1] \to (c,d) \in \mathcal{R}$ is twice differentiable with

$$\int_0^1 \{m^{(2)}(x)\}^2\,dx < \infty \quad \text{and} \quad m^{(p)}(0) = m^{(p)}(1), \qquad p = 0,1,2,\ldots,$$

where $m^{(p)}$ denotes the $p$th derivative of $m$.
(A4) The bandwidth sequence $h$ depends on $n$ and satisfies $h \to 0, nh \to \infty$ as $n \to \infty$.
(A5) The function $\rho$ is a continuous function symmetric about zero and is differentiable everywhere except possibly at a finite number of points.
(A6) There are constants $c_0, c_1 > 0$ such that, for $x \in [0,1]$,

$$|E_F \psi(y - m(x) + s)| > c_0|s|, |s| < c_1,$$

where $\psi$ is the derivative of $\rho$.

(A7) There exists a constant $c_2$ such that, for $x \in [0,1]$,
$$|E_F \psi(y - m(x) + s)| \leq c_2 |s|.$$

(A8) $E_F[\psi(\cdot)]^2 < \infty; E_F[\psi'(\cdot)]^2 < \infty.$

REMARK. For a symmetric noise distribution $F$, assumptions (A6) and (A7) are satisfied for $\psi(u) = u$. For nonlinear $\psi$ functions, (A6) is satisfied if $E_F \psi'(y - m(x) + s) > c_0$ and (A7) is satisfied if $E_F \psi'(y - m(x) + s) < c_2$ for small $s$. Therefore, these conditions ensure that $E_F \psi'(y - m(x) + s)$ is a positive bounded value.

Hereafter, the notation is simplified by writing $m_i$ instead of $m(x_i)$; $\tilde{m}_i$ instead of $\tilde{m}(x_i)$; $\tilde{m}_{-i}$ instead of $\tilde{m}_{-i}(x_i)$; and $E$ instead of $E_F$. Furthermore, the names $c$ and $q$ are generic names for constants; they may represent different values in different contexts.

One difficulty in working with the robust estimate, $\tilde{m}_i$, is the representation of it in a workable form. From the definition of $\tilde{m}_i$, it can expanded in a Taylor series about $m_i$ as

$$(6) \quad \sum_{j=1}^n \alpha_j(x_i) \tilde{\psi}(y_j - m_i) + \sum_{j=1}^n \alpha_j(x_i)(m_i - \tilde{m}_i) \tilde{\psi}'(y_j - m_i + \nu_i) = 0,$$

where $|\nu_i| < |m_i - \tilde{m}_i|$. Therefore, if there exists a constant $q > 0$ such that

$$(7) \quad \sup_h \left| \sum_{j=1}^n \alpha_j(x_i) \tilde{\psi}'(y_j - m_i + \nu_i) - q \right| \to 0 \quad \text{a.s.},$$

then, on the event $\mathcal{P} = \{\inf_h \sum_{j=1}^n \alpha_j(x_i) \tilde{\psi}'(y_j - m_i + \nu_i) > 0\}$, which occurs a.s., $\tilde{m}_i$ can be represented as

$$(8) \quad m_i - \frac{\sum_{j=1}^n \alpha_j(x_i) \tilde{\psi}(y_j - m_i)}{q},$$

which is a linear combination of the random variables $\tilde{\psi}(y_j - m_i)$. This result can be shown by using (A6), (A7) and making use of Theorem 2 in [13] and Theorem 2 in [28].

Similarly, $\tilde{m}_{-i}$ can be represented as

$$(9) \quad m_i - \frac{\sum_{j \neq i}^n \alpha_j(x_i) \tilde{\psi}(y_j - m_i)}{q}.$$

PROPOSITION 1. *If the conditions* (A1)–(A8) *are satisfied, then on the event* $\mathcal{P}$ *and over an interval* $[\zeta_1 n^{-1/5}, \zeta_2 n^{-1/5}]$ *for some suitable* $\zeta_1, \zeta_2$,

$$(10) \quad RCRVD(h) = c + \frac{E\psi'(\cdot)}{2} MASE(h) + o_p(n^{-4/5})$$

*uniformly in* $h$, *where $c$ is a constant w.r.t. $h$.*

Proposition 1 shows that RCRVD differs from MASE by a constant shift and a constant multiple. Neither of these is dependent on the bandwidth. Hence, asymptotically, minimizing RCRVD and MASE w.r.t. $h$ are equivalent. This result is shown in Proposition 3 below.

PROPOSITION 2. *If the conditions* (A1)–(A8) *are satisfied, then*
$$ASE(h) = MASE(h) + o_p(n^{-4/5})$$
*uniformly in $h$.*

PROPOSITION 3. *If the conditions* (A1)–(A8) *are satisfied,*
$$h_{\text{CRVD}} = \arg\min_h RCRVD(h),$$
$$h_{\text{ASE}} = \arg\min_h ASE(h),$$
$$h_{\text{MASE}} = \arg\min_h MASE(h)$$
*are all equivalent as $n \to \infty$. In particular, they are all equal to $kn^{-1/5}$, where*

(11) $$k = \left[\frac{\int K^2(u)\,du E\tilde{\psi}^2(\cdot)}{\int_0^1 \{m^{(2)}(x)\}^2\,dx \int u^2 K(u)\,du \{E\tilde{\psi}'(\cdot)\}^2}\right]^{1/5}.$$

Note that $k$ depends on $\tilde{\psi}$ rather than $\psi$, and hence the bandwidth selected by RCRVD is asymptotically independent of the choice of $\psi$ (or $\rho$).

**3. Simulation study results.** A simulation study on the finite sample performance of RCRVD was performed. In this study the following regression function, $m(x)$, was used:

(12) $$m(x) = \sin(2\pi x), \qquad 0 < x < 1.$$

The observations $Y_i$ were taken at $x_i = i/n$, for $n = 257$, and $\varepsilon_i$ was an error from one of the following three distributions: (1) $N(0, 0.2)$; (2) Contaminated normal $0.9N(0, 0.2) + 0.1N(0, 1.8)$ and (3) Contaminated normal $0.8N(0, 0.2) + 0.2N(0, 3)$, where $N(\mu, \sigma)$ stands for the normal distribution. The program for computing the robust smoother and the cross-validation rules was written in Fortran 90, incorporating the algorithm of Härdle [8]. Since Härdle's algorithm used a fast Fourier transform in the computations, the function was evaluated on an interval with the number of grids a power of 2, which was set to an upper limit of 1024. Hence, the choice of the sample size of 257 ($= 2^8 + 1$) in the simulation was for computational convenience. Of course, the method works for any other sample size. A standard normal kernel was used for all smoothing. Robustness in smoothing was achieved by setting $\tilde{\rho}$ to Huber's $\rho$ [16],

$$\rho(u) = \begin{cases} \dfrac{u^2}{2}, & \text{if } |u| < c, \\ c|u| - \dfrac{c^2}{2}, & \text{if } |u| \geq c, \end{cases}$$

where the threshold $c = 0.5$ was chosen to obtain the appropriate degree of robustness.

When using RCRVD, a second $\rho$, not necessarily the same as $\tilde{\rho}$, is used to discount the extreme prediction errors. In this study three different choices of $\rho$ were considered: (1) Huber's $\rho$ with $c = 0.5$; (2) Huber's $\rho$ with $c = 1$; (3) $\rho(u) = |u|$. The choice $\rho(u) = |u|$ is equivalent to the $L_1$ cross-validation method of Wang and Scott [26]. We note that, asymptotically, the different choices of $\rho$ all give consistent bandwidths. But in finite samples, using different choices of $\rho$ may lead to different bandwidths.

First, the behaviors of the RCRVDs and CRVD were studied in a simulated example. One sample of 257 observations was generated using each of the three error distributions described above. The data, along with the robust smoother ($\tilde{m}$) using the optimal bandwidth, $h_{\text{ASE}}$, were plotted in Figures 1(a), 2(a) and 3(a). Additionally, the values of the RCRVDs, CRVD and ASE were plotted as functions of the bandwidth, $h$, in Figures 1(b), 2(b) and 3(b). Since the optimal bandwdith chosen by a particular bandwidth selector is not affected by a constant shift in the value of the selector, in Figures 1(b), 2(b) and 3(b) the plots of CRVD and RCRVD have been shifted by constant amounts for ease of comparison. When there were no outliers [Figure 1(b)], the minima of all bandwidth selectors (CRVD and all three RCRVDs) were close to that of ASE. When the data were contaminated [Figures 2(b) and 3(b)], the minima of CRVD were very different from those using ASE. The minima of the RCRVDs, on the other hand, were still very similar to those of ASE, even under heavy contamination [Figure 3(b)].

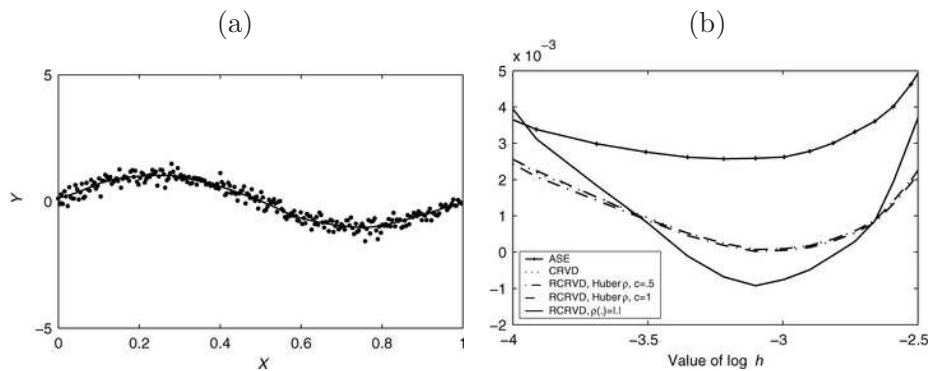

Fig. 1. (a) $N(0, 0.2)$ data with robust smoother using $h_{\text{ASE}}$. (b) Plot of ASE, CRVD and RCRVDs vs. $\log h$.

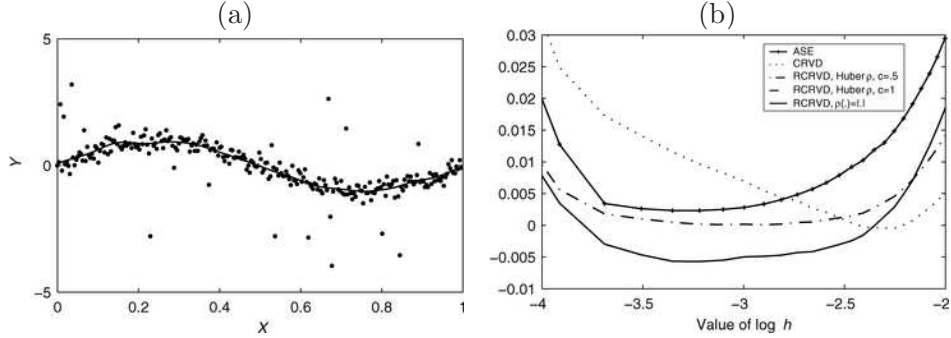

FIG. 2. (a) $0.9N(0,0.2)+0.1N(0,1.8)$ *data with robust smoother using* $h_{\text{ASE}}$. (b) *Plot of ASE, CRVD and RCRVDs vs.* $\log h$.

In the simulation study 100 samples of 257 observations were generated from each of the three noise distributions described above. Eight different bandwidth selection methods were considered in this study. These methods included the standard CRVD and the three different versions of RCRVD considered above. In addition, four selectors based on the robust plug-in method of Boente, Fraiman and Meloche [2] were included. The plug-in method they considered was based on finding $h_{\text{PLUG-IN}} = kn^{-1/5}$, where $k$ is given by (11). Following Boente, Fraiman and Meloche [2], in the expression for $h_{\text{PLUG-IN}}$, $E\tilde{\psi}^2(\cdot)/E^2\tilde{\psi}'(\cdot)$ was replaced in (11) by

$$\hat{\sigma}_R = \left(\frac{1}{0.6745\sqrt{2}}\operatorname{med}_i |Y_{i+1} - Y_i|\right)^2,$$

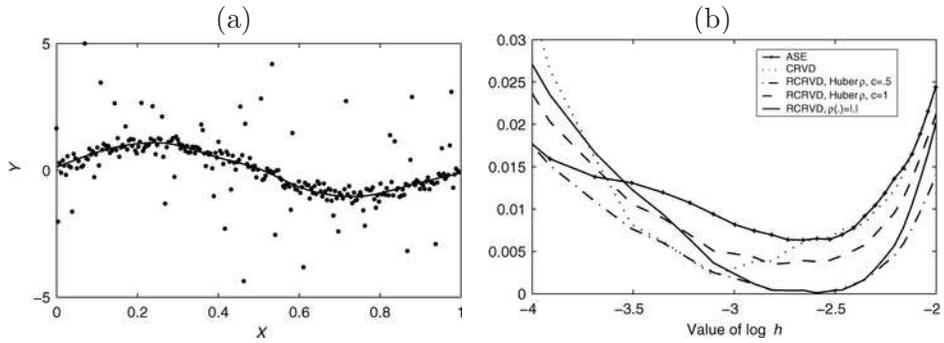

FIG. 3. (a) $0.8N(0,0.2)+0.2N(0,3)$ *data with robust smoother using* $h_{\text{ASE}}$. (b) *Plot of ASE, CRVD and RCRVDs vs.* $\log h$.

where med$_i$ stands for median over $i$, $i = 1, \ldots, n - 1$. The quantity $\int_0^1 \{m^{(2)}(x)\}^2 \, dx$ was evaluated by estimating $m^{(2)}(x)$ on a fine grid and then applying numerical integration. The quantity $m^{(2)}(x)$ was estimated using a robust kernel regression for derivatives ([2], page 119) using a "pilot bandwidth," $h_0$, which had to be determined subjectively. To examine the sensitivity of the plug-in method to the choice of $h_0$, four choices of $h_0$, $0.02, 0.03, 0.04$ and $0.06$ were considered.

For each sample in the simulation, we calculated the quantity

$$\text{(13)} \qquad \frac{|h_{\text{method}} - h_{\text{ASE}}|}{h_{\text{ASE}}},$$

where method = CRVD, RCRVD or PLUG-IN. Ideally, a method should choose a bandwidth close to that chosen by ASE. Therefore, for a good method, the distribution of (13) should have most of its density clustered around zero. Based on the results from the simulations, the frequency distribution of (13) was tabulated for the different methods (Table 1). Each entry in Table 1 gives the number of times out of 100 simulations that the quantity (13) fell within the interval given in the column heading. For example, for CRVD, under $N(0, 0.2)$ data, there were 33 times out of 100 simulations that (13) was smaller than 0.1.

The results of the simulations showed that when the data were normal [Table 1(a)], all cross-validation based methods (CRVD and RCRVD) behaved similarly and satisfactorily, that is, the quantity (13) was less than 0.6 in most of the simulations, indicating that the methods chose bandwidths close to those given by ASE in most of the simulation runs. The performances of the plug-in methods, however, were quite different. Among the four plug-in methods, only the one using a pilot bandwidth, $h_0 = 0.04$, gave comparable performance to the cross-validation methods, while the other three gave very disappointing results. This pattern of performance by the plug-in method was consistent throughout the study, giving rise to implications later.

For mildly contaminated normal errors [Table 1(b)], the cross-validation methods started to diverge in their performances. CRVD often chose a bandwidth that was much smaller than that selected by ASE, resulting in moderately large values of (13) in a lot of the samples. For RCRVD, there was no evidence of any impaired performance using either Huber's $\rho$ with $c = 0.5$ or $\rho(\cdot) = |\cdot|$, even though the performance using Huber's $\rho$ with $c = 1$ yielded less satisfactory results. The performance of CRVD became worse when the proportion and severity of contamination were high [Table 1(c)]. Among the RCRVDs, the performance was better in the selector with $\rho(u) = |u|$.

To gain more insight into the behaviors of the different methods, the Monte Carlo mean and standard deviation of the optimal bandwidth chosen by the different methods are summarized in Table 2. From this table the following observations can be noted. First, despite the poor performance of

TABLE 1
*Frequency distribution of $|h_{\mathrm{method}} - h_{\mathrm{ASE}}|/h_{\mathrm{ASE}}$ for the various bandwidth selectors: 100 data sets of size 257 with three different error distributions*

| Method | Value of $|h_{\mathrm{method}} - h_{\mathrm{ASE}}|/h_{\mathrm{ASE}}$ | | | | | | | |
|---|---|---|---|---|---|---|---|---|
| | ≤0.1 | 0.1–0.2 | 0.2–0.4 | 0.4–0.6 | 0.6–0.8 | 0.8–1 | 1–2 | >2 |
| (a) Error distribution: $N(0, 0.2)$ | | | | | | | | |
| CRVD | 33 | 20 | 33 | 8 | 6 | 0 | 0 | 0 |
| RCRVD, Huber's $\rho$, $c = 0.5$ | 32 | 21 | 33 | 6 | 8 | 0 | 0 | 0 |
| RCRVD, Huber's $\rho$, $c = 1$ | 33 | 20 | 33 | 8 | 6 | 0 | 0 | 0 |
| RCRVD, $\rho(u) = |u|$ | 27 | 27 | 31 | 12 | 2 | 0 | 1 | 0 |
| PLUG-IN, $h_0 = 0.02$ | 0 | 0 | 14 | 68 | 18 | 0 | 0 | 0 |
| PLUG-IN, $h_0 = 0.03$ | 18 | 18 | 49 | 15 | 0 | 0 | 0 | 0 |
| PLUG-IN, $h_0 = 0.04$ | 38 | 31 | 27 | 4 | 0 | 0 | 0 | 0 |
| PLUG-IN, $h_0 = 0.06$ | 17 | 18 | 26 | 19 | 11 | 5 | 4 | 0 |
| (b) Error distribution: $0.9N(0, 0.2) + 0.1N(0, 1.8)$ | | | | | | | | |
| CRVD | 20 | 13 | 28 | 19 | 8 | 5 | 7 | 0 |
| RCRVD, Huber's $\rho$, $c = 0.5$ | 43 | 17 | 25 | 11 | 4 | 0 | 0 | 0 |
| RCRVD, Huber's $\rho$, $c = 1$ | 27 | 26 | 27 | 10 | 6 | 4 | 0 | 0 |
| RCRVD, $\rho(u) = |u|$ | 36 | 25 | 23 | 11 | 3 | 2 | 0 | 0 |
| PLUG-IN, $h_0 = 0.02$ | 0 | 0 | 3 | 40 | 56 | 1 | 0 | 0 |
| PLUG-IN, $h_0 = 0.03$ | 6 | 7 | 36 | 51 | 0 | 0 | 0 | 0 |
| PLUG-IN, $h_0 = 0.04$ | 23 | 25 | 43 | 9 | 0 | 0 | 0 | 0 |
| PLUG-IN, $h_0 = 0.06$ | 28 | 26 | 29 | 11 | 6 | 0 | 0 | 0 |
| (c) Error distribution: $0.8N(0, 0.2) + 0.2N(0, 3)$ | | | | | | | | |
| CRVD | 5 | 15 | 18 | 25 | 17 | 10 | 8 | 2 |
| RCRVD, Huber's $\rho$, $c = 0.5$ | 36 | 26 | 21 | 5 | 8 | 1 | 3 | 0 |
| RCRVD, Huber's $\rho$, $c = 1$ | 18 | 20 | 33 | 16 | 8 | 2 | 3 | 0 |
| RCRVD, $\rho(u) = |u|$ | 38 | 29 | 17 | 6 | 6 | 0 | 4 | 0 |
| PLUG-IN, $h_0 = 0.02$ | 0 | 0 | 1 | 12 | 71 | 16 | 0 | 0 |
| PLUG-IN, $h_0 = 0.03$ | 0 | 4 | 14 | 52 | 29 | 1 | 0 | 0 |
| PLUG-IN, $h_0 = 0.04$ | 8 | 7 | 31 | 50 | 3 | 1 | 0 | 0 |
| PLUG-IN, $h_0 = 0.06$ | 23 | 32 | 32 | 9 | 2 | 0 | 2 | 0 |

CRVD, its means were very similar to the corresponding means under ASE. However, its standard deviation, especially under heavy contamination (last column of Table 2), was much higher than those of ASE and the other selectors. This result indicates that CRVD sometimes chose bandwidths that were quite far away from those using ASE. Second, the standard deviations of all plug-in methods were small, which is consistent with results seen in other applications of the plug-in method (cf. [6]). However, for these methods the means were very sensitive to the choice of the pilot bandwidth, $h_0$. These observations indicate that the poor performance of the plug-in methods results from the bias that can arise from an inappropriate choice of $h_0$. Therefore, for the plug-in methods a good choice of pilot bandwidth is cruical. Third, in contrast to the other methods, all the RCRVDs seemed to choose unbiased bandwidths, and the standard deviations, even though higher than those of the plug-in methods, were not unreasonably large.

TABLE 2
*Mean and standard deviation of the optimal bandwidths for the various bandwidth selectors:*
100 *data sets of size* 257, *with three different error distributions*

| | Error distribution | | |
|---|---|---|---|
| | $N(0, 0.2)$ | $0.9N(0, 0.2)$ $+ 0.1N(0, 1.8)$ | $0.8N(0, 0.2)$ $+ 0.2N(0, 3)$ |
| **Method** | Mean (SD) | Mean (SD) | Mean (SD) |
| ASE | 0.0471 (0.0000507) | 0.0575 (0.0001232) | 0.0735 (0.0002968) |
| CRVD | 0.0488 (0.0000457) | 0.0593 (0.0003014) | 0.0675 (0.0026117) |
| RCRVD, Huber's $\rho$, $c = 0.5$ | 0.0481 (0.0000558) | 0.0581 (0.0000799) | 0.0744 (0.0001589) |
| RCRVD, Huber's $\rho$, $c = 1$ | 0.0488 (0.0000460) | 0.0585 (0.0000976) | 0.0703 (0.0002314) |
| RCRVD, $\rho(u) = |u|$ | 0.0477 (0.0000847) | 0.0588 (0.0001253) | 0.0758 (0.0001811) |
| PLUG-IN, $h_0 = 0.02$ | 0.0228 (0.0000058) | 0.0222 (0.0000084) | 0.0199 (0.0000199) |
| PLUG-IN, $h_0 = 0.03$ | 0.0349 (0.0000101) | 0.0346 (0.0000090) | 0.0327 (0.0000245) |
| PLUG-IN, $h_0 = 0.04$ | 0.0461 (0.0000175) | 0.0466 (0.0000185) | 0.0437 (0.0000559) |
| PLUG-IN, $h_0 = 0.06$ | 0.0620 (0.0000478) | 0.0659 (0.0000448) | 0.0665 (0.0000771) |

**4. Application to water-quality data.** The tremendous development in Florida from the 1960s to the 1990s and the changes in the demand in and practice of water use over that period had caused a decrease in water supply and quality in the 1990s. In response to these problems, in 1996 the U.S. Geological Survey conducted a study on the long term water quality trend in Southern Florida. In this section the methods considered in this paper are applied to data from the U.S. Geological Survey. The data were collected at two discharge stations—one within Big Cypress National Preserve and one near Biscayne Bay [18].

A large number of water-quality constituents and flow data were collected periodically between 1966–1994 and used in the survey. The methods in this paper are illustrated using the dissolved solids data collected at Tamiami Canal station inside the Big Cypress Preserve. The data included the data that were used in the survey and also data collected up to 1999. The total number of observations was 118. The raw data of dissolved solids level (in mg/l), as a function of time, were plotted in Figure 4(a). The plot clearly shows outliers in the data. A robust smoother was fit by setting $\tilde{\rho}$ to Huber's $\rho$. The value of $c$ in $\tilde{\rho}$ was obtained by first fitting a LOWESS curve [3] to the data. A smoothing parameter of 0.5 was used for the LOWESS fit, which was the value used by Lietz [18]. The value of $c = 1.345 \times (\text{med}_i |\text{res}_i|)/0.645$, where $\text{res}_i$ was the residual of the $i$th observation from the LOWESS plot, was used for $\tilde{\rho}$. Using the data, the value of $c = 64$ was obtained.

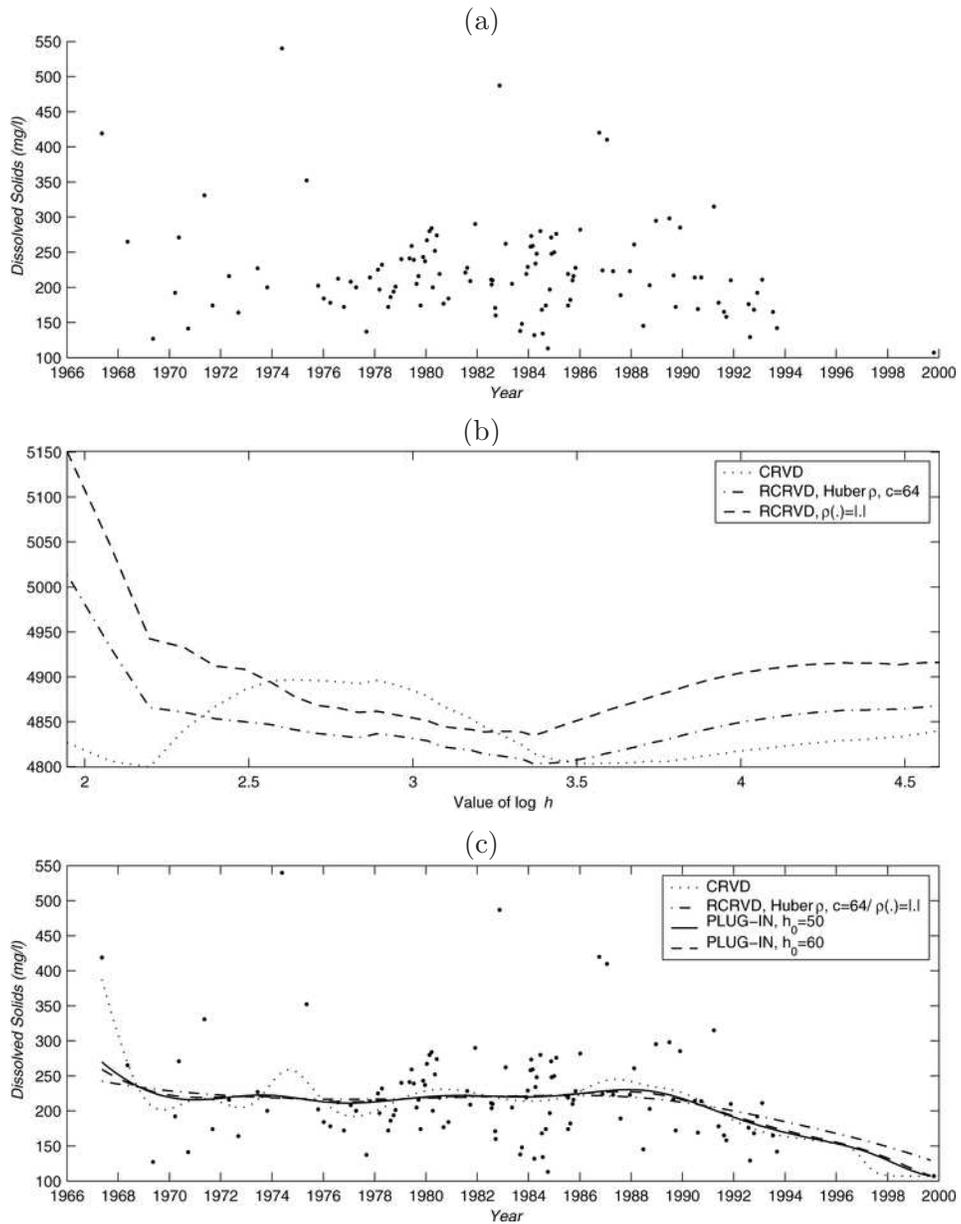

Fig. 4. (a) *Scatterplot of dissolved solids as a function of time at the Tamiami Canal station.* (b) *Plot of CRVD and RCRVDs vs.* $\log h$. (c) *Smoothers for dissolved solids as a function of time at the Tamiami Canal station.*

The following methods for bandwidth selection were considered: (1) CRVD; (2) RCRVD with Huber's $\rho, c = 64$; (3) RCRVD with $\rho(u) = |u|$; (4) the plug-in method with pilot bandwidth $h_0 = 50$; and (5) the plug-in method with $h_0 = 60$.

The robust smoother was evaluated on a grid with intervals of 50 days. Therefore, all results were quoted in units of 50 days. Figure 4(b) gives the CRVD and the two RCRVDs as functions of the bandwidth. It is clear that the minimum of CRVD appeared much earlier than the other two curves. The optimal bandwidths, obtained by grid searches, were 9 ($= 9 \times 50$ days), 29 and 29, respectively, for CRVD, RCRVD with Huber's $\rho, c = 64$, and RCRVD with $\rho(u) = |u|$. The smaller minimum of CRVD made sense; the few extreme values in the dataset caused CRVD to suggest a rougher curve to downweight the extreme prediction errors. The optimal plug-in bandwidths were 18 and 20.4, respectively, for $h_0 = 50$ and 60.

The robust smoothers based on the different bandwidths are plotted in Figure 4(c). The rough curve (dotted line) is based on the bandwidth using CRVD. The other curves are similar to each other. The results show a downward trend in dissolved solids over time at Tamiami Canal station.

**5. Conclusion.** This study demonstrates that, in data suspected of containing outliers, using standard cross-validation can lead to bandwidths that are not optimal with respect to the usual error criteria. On the other hand, the robust cross-validation method suggested in this article provides bandwidths that are asymptotically optimal with respect to these criteria. Proposition 1 demonstrates that, asymptotically, RCRVDs using different forms of $\rho(\cdot)$ differ only by a constant. This finding suggests that *all* robust cross-validation methods that satisfy the assumptions of the proposition give the same asymptotic bandwidth. Asymptotically, the bandwidth chosen is optimal with respect to MASE and ASE. The assumptions on the form of $\rho(\cdot)$ are all very mild. The most important assumption is the symmetry of $\rho(\cdot)$, which ensures deviations from either side of the unknown curve to be weighted equally. In small and moderate samples, however, different forms of $\rho(\cdot)$ will give different results. From the derivation of Proposition 1, it is important to note that the extent of how well RCRVD approximates MASE or ASE depends on how well $\tilde{m}_i$ and $\tilde{m}_{-i}$ are represented as (8) and (9), respectively. It has been shown that a $\rho$ with a small value of $\int \{\psi^{(2)}(\cdot)\}^2$ will give satisfactory results.

The method suggested here includes the $L_1$ robust bandwidth selector of Wang and Scott [26] as a special case. Another alternative to the proposed method is the plug-in method suggested by Boente, Fraiman and Meloche [2]. The simulations in this study show that, the bandwidth chosen using a plug-in method, despite having a smaller standard deviation than the cross-validation method, is highly sensitive to the choice of the pilot bandwidth. Similar observations were made by Boente, Fraiman and Meloche [2], Tables 1–4, though the sensitivity they observed was less severe.

The method described here assumes a circular design in which there are no boundary problems. But the method can be easily adapted to work for a general regression function with the incorporation of a weight function that discounts the contribution of the boundary observations to the cross-validation; see, for example, [15]. The results for Propositions 1–3 will still hold under this situation.

The method studied in this article is a general robust method. It can be applied to problems other than kernel smoothing where cross-validation is needed. For example, Ronchetti, Field and Blanchard [24] considered robust cross-validation for model selection in a linear regression. This method can be used in that application.

## APPENDIX: SKETCHES OF THE PROOFS

PROOF OF PROPOSITION 1. Without loss of generality, it is assumed here that $\tilde{\rho} = \rho$. This assumption will simplify the notation in the proof considerably. The result of the proposition still holds when $\tilde{\rho} \neq \rho$. The first step of the proof is to expand RCRVD in Taylor series. To show (10), it is required to demonstrate that the sup norms of the l.h.s. of (A.1)–(A.3) over the interval $(\zeta_1 n^{-1/5}, \zeta_2 n^{-1/5})$ vanish. The general idea is to partition $(\zeta_1 n^{-1/5}, \zeta_2 n^{-1/5})$ into many small subintervals and then "discretize" the sup norm problem into one of finding the maximum of a finite number of sup norms within these subintervals, which themselves can be bounded easily by established results [see, e.g., (A.5)–(A.8)],

$$RCRVD(h) = n^{-1} \sum_{i=1}^{n} \rho(\varepsilon_i + m_i - \tilde{m}_{-i});$$

expand $\rho$ by Taylor series at the points $\varepsilon_i, i = 1, \ldots, n$,

$$RCRVD(h) = n^{-1} \sum_{i=1}^{n} \rho(\varepsilon_i) + D_1(h) + D_2(h) + D_3(h),$$

where

$$D_1(h) = n^{-1} \sum_{i=1}^{n} (m_i - \tilde{m}_{-i}) \psi(\varepsilon_i),$$

$$D_2(h) = n^{-1} \sum_{i=1}^{n} \frac{(m_i - \tilde{m}_{-i})^2}{2} \psi'(\varepsilon_i)$$

and $D_3(h) = n^{-1}\sum_{i=1}^n R_i$ are the remainder terms in the Taylor series expansions. To prove (10), it is required to demonstrate that

(A.1) $$D_1(h) = o_p(n^{-4/5}),$$

(A.2) $$D_2(h) - \frac{E\psi'(\cdot)}{2} MASE(h) = o_p(n^{-4/5}),$$

(A.3) $$D_3(h) = o_p(n^{-4/5})$$

uniformly in $h$.

$D_1(h)$ can be written as

(A.4) $$D_1(h) = q^{-1}n^{-1}\sum_{i=1}^n \sum_{j\neq i} \alpha_j(x_i)\psi(\varepsilon_j)\psi(\varepsilon_i).$$

Therefore, to show (A.1), for some $\varepsilon > 0$,

(A.5) $$P\left(\sup_h |D_1(h)| \geq n^{-4/5}\varepsilon\right)$$
$$\leq P\left(\sup_{e=2,\ldots,m} \sup_{s_{e-1} < h_e < s_e} |D_1(h_e) - D_1(s_e)| \geq n^{-4/5}\varepsilon/2\right)$$
$$+ P\left(\sup_{e=2,\ldots,m} |D_1(h_e)| \geq n^{-4/5}\varepsilon/2\right),$$

where the intervals $s_1 \leq s_2 \leq \cdots \leq s_m$ form a partition of the interval $(\zeta_1 n^{-1/5}, \zeta_2 n^{-1/5})$ and the sizes of the intervals $(s_{e-1}, s_e), e = 2, \ldots, m$, are to be chosen small enough. In other words, if $m$ is large enough, the first term on the r.h.s. of (A.5) becomes negligible compared to the second term. Writing $\eta = \varepsilon/2$, it is necessary to demonstrate

(A.6) $$P\left(\sup_{e=2,\ldots,m} |D_1(h_e)| \geq n^{-4/5}\eta\right)$$

goes to zero. By Bonferroni's inequality, (A.6) is bounded above by

(A.7) $$m \sup_{e=2,\ldots,m} P(|D_1(h_e)| \geq n^{-4/5}\eta) \leq m \sup_{e=2,\ldots,m} E(n^{4/5}\eta^{-1}|D_1(h_e)|)^{2k},$$

where $k = 1, 2, \ldots$.

Define $z_i = \psi(\varepsilon_i), i = 1, \ldots, n$. $D_1(h_e)$ is a quadratic form of the $z_i$'s, $\beta_{ij} z_i z_j$, where $\beta_{ij} = q^{-1}n^{-1}\alpha_j(x_i), \beta_{ii} = 0$ and $z_i, z_j$ are independent for $i \neq j$. Now, using Theorem 2 of [28],

$$E(n^{4/5}\eta^{-1}|D_1(h_e)|)^{2k}$$
$$\leq c\bigg[n^{8/5}(\eta q n h_e)^{-2}$$

$$\text{(A.8)} \quad \times \sum_{i=1}^{n}\sum_{j\neq i}\left\{(nh_e)^{-1}K\left(\frac{x_i - x_j}{h_e}\right)\bigg/\sum_{l=1}^{n}K\left(\frac{x_i - x_l}{h_e}\right)\right\}^2$$

$$\times E\psi^2(\varepsilon_i)E\psi^2(\varepsilon_j)\Bigg]^k$$

$$\leq c\left[n^{(-4+8/5)k}h_e^{-2k}\sum_{l=2}^{2k}n^l h_e^l\{E\psi^2(\varepsilon_i)\}^{2k}\right].$$

Note that, in the above the Nadaraya–Watson form of $\alpha_i(x)$ is used. But the result still holds for Gasser–Müller, Priestley–Chao or other weights. Therefore,

$$m\sup_{e=2,\ldots,m} E(n^{4/5}\eta^{-1}|D_1(h_e)|)^{2k}$$

$$\leq m\sup_{e=2,\ldots,m} c\left[n^{(-4+8/5)k}h_e^{-2k}\sum_{l=2}^{2k}n^l h_e^l\{E\psi^2(\varepsilon_i)\}^{2k}\right],$$

which goes to zero for fixed $m$, $h_e \in [s_{e-1}, s_e]$ and $n$ going to $\infty$, if $k$ is chosen to be sufficiently large.

The derivation of (A.2) is similar to that for (A.1). First write

$$D_2(h) - ED_2(h) = \frac{1}{2n}\sum_{i=1}^{n}(m_i - \tilde{m}_i)^2\psi'(\varepsilon_i) - \sum_{i=1}^{n}E(m_i - \tilde{m}_{-i})^2 E\psi'(\varepsilon_i)$$

$$= \frac{1}{2n}\sum_{i=1}^{n}\{(m_i - \tilde{m}_i)^2 - E(m_i - \tilde{m}_{-i})^2\}\psi'(\varepsilon_i)$$

$$- \sum_{i=1}^{n}E(m_i - \tilde{m}_{-i})^2\{\psi'(\varepsilon_i) - E\psi'(\varepsilon_i)\}$$

$$= \frac{1}{2}\left[n^{-1}\sum_{i=1}^{n}W_1^i(h) + n^{-1}\sum_{i=1}^{n}W_2^i(h)\right].$$

It will be demonstarted that $n^{-1}\sum_{i=1}^{n}W_1^i = o_p(n^{-4/5})$ and $n^{-1}\sum_{i=1}^{n}W_2^i = o_p(n^{-4/5})$ uniformly in $h$. Since $\psi(\cdot)$ is bounded, for $n^{-1}\sum_{i=1}^{n}W_1^i$ it is enough to show that $[(m_i - \tilde{m}_i)^2 - E(m_i - \tilde{m}_{-i})^2]$, $i = 1, \ldots, n$, are $o_p(n^{-4/5})$. Each term is expanded as

$$\left[\frac{\sum_{j\neq i}\alpha_j(x_i)\psi(y_j - m_i)}{q}\right]^2 - E\left[\frac{\sum_{j\neq i}\alpha_j(x_i)\psi(y_j - m_i)}{q}\right]^2$$

$$= \left\{\sum_{j\neq i}\sum_{k\neq i}\alpha_j(x_i)\alpha_k(x_i)\{\psi(y_j - m_i) - E\psi(y_j - m_i)\}\right.$$

$$\times \{\psi(y_k - m_i) - E\psi(y_k - m_i)\}\bigg\}\bigg/q^2$$

$$+ \frac{\sum_{j\neq i}\sum_{k\neq i}\alpha_j(x_i)\alpha_k(x_i)\psi(y_j - m_i)E\psi(y_k - m_i)}{q^2}$$

(A.9) $$+ \frac{\sum_{j\neq i}\sum_{k\neq i}\alpha_j(x_i)\alpha_k(x_i)\psi(y_k - m_i)E\psi(y_j - m_i)}{q^2}$$

$$- \frac{\sum_{j\neq i}\alpha_j^2(x_i)[E\psi^2(y_j - m_i)\{E\psi(y_k - m_i)\}^2]}{q^2}$$

$$- \frac{2\sum_{j\neq i}\sum_{k\neq i, k\neq j}\alpha_j(x_i)\alpha_k(x_i)E\psi(y_j - m_i)E\psi(y_k - m_i)}{q^2}$$

$$= \bigg\{\sum_{j\neq i}\sum_{k\neq i}\alpha_j(x_i)\alpha_k(x_i)\{\psi(y_j - m_i) - E\psi(y_j - m_i)\}$$

$$\times \{\psi(y_k - m_i) - E\psi(y_k - m_i)\}\bigg\}\bigg/q^2$$

$$+ \frac{\sum_{j\neq i}\sum_{k\neq i}\alpha_j(x_i)\alpha_k(x_i)\{\psi(y_j - m_i) - E\psi(y_j - m_i)\}E\psi(y_k - m_i)}{q^2}$$

$$+ \frac{\sum_{j\neq i}\sum_{k\neq i}\alpha_j(x_i)\alpha_k(x_i)\{\psi(y_k - m_i) - E\psi(y_k - m_i)\}E\psi(y_j - m_i)}{q^2}$$

$$- \frac{\sum_{j\neq i}\alpha_j(x_i)\alpha_k(x_i)[E\psi^2(y_j - m_i)\{E\psi(y_k - m_i)\}^2]}{q^2}$$

$$= W_{11}^i(h) + W_{12}^i(h) + W_{13}^i(h) + W_{14}^i(h).$$

To bound $W_{11}^i(h)$, write $r_j = \psi(y_j - m_i) - E\psi(y_j - m_i)$. Note that $r_i, r_j, i \neq j$, are independent r.v.'s with zero mean. An application of Theorem 2 of [28] gives

(A.10) $$P\bigg(\sup_h |W_{11}^i(h)| \geq n^{-4/5}\varepsilon\bigg) \leq cmn^{(-4+8/5)k}h_e^{-2k}\sum_{l=2}^k (nh_e)^l.$$

Similarly, $W_{12}^i(h), W_{13}^i(h), W_{14}^i(h)$ can also be bounded as in (A.9). Therefore, $W_{11}^i(h), W_{12}^i(h), W_{13}^i(h), W_{14}^i(h)$ all vanish as long as a large enough $k$ is chosen, as $n \to \infty$ and $h \in (\zeta_1 n^{-1/5}, \zeta_2 n^{-1/5})$, with $m$ fixed. From [10], $|E(m_i - \tilde{m}_{-i})^2 - E(m_i - \tilde{m}_i)^2| = o_p(n^{-4/5})$ and $E(m_i - \tilde{m}_i)^2 = O_p(n^{-1}h^{-1} + h^4)$. Therefore, for $\varepsilon > 0$,

$$P\bigg(\sup_h \bigg|n^{-1}\sum_{i=1}^n W_2^i(h)\bigg| \geq n^{-4/5}\varepsilon\bigg)$$

(A.11)
$$\leq c\left[m(n^{-1}h_e^{-1} + h_e^4)^k n^{(-4+8/5)k} \sum_{l=2}^{k}(nh_e)^l\right],$$

which goes to zero as $k$ is chosen large enough. Hence, $W_2^i(h) = o_p(n^{-4/5})$, uniformly in $h$. Now, since $E(D_2) = E\psi'(\cdot)MASE(h)/2 + o_p(n^{-4/5})$, (A.2) follows from the above derivations. Finally, (A.3) can be easily shown by the same techniques that are used to show (A.1) and (A.2). □

PROOF OF PROPOSITION 2. The proof follows the same route as that of Theorem 1 of [19]. □

PROOF OF PROPOSITION 3. By standard technique, it can be shown that

(A.12)
$$MASE(h) = (nh)^{-1}\int K^2(u)\,du \frac{E\tilde{\psi}^2(\cdot)}{\{E\tilde{\psi}'(\cdot)\}^2}$$
$$+ \frac{h^4}{4}\int_0^1 [m^{(2)}(x)]^2\,dx\left(\int u^2 K(u)\,du\right)^2$$
$$+ o(n^{-1}h^{-1} + h^4)$$

and $\arg\min_h MASE(h)$ is given by $kn^{-1/5}$, where $k$ is as in (11). Using a similar approach as in [23],
$$D(\varepsilon) = \inf_{|u-v|>n^{-1/5}\varepsilon} n^{-4/5}|RCRVD(u) - RCRVD(v)|.$$

It follows that

(A.13)
$$P(|h_{\text{CRDV}} - h_{\text{MASE}}| > n^{-1/5}\varepsilon)$$
$$\leq P\left(n^{4/5}\sup_h |RCRVD(h_{\text{MASE}}) - RCRVD(h_{\text{RCRVD}})| > D(\varepsilon)\right)$$
$$\leq P\left(n^{4/5}\sup_h \left|RCRVD(h_{\text{MASE}}) - \frac{E\psi'(\cdot)}{2}MASE(h_{\text{MASE}})\right.\right.$$
$$\left.\left. + \frac{E\psi'(\cdot)}{2}MASE(h_{\text{MASE}}) - RCRVD(h_{\text{RCRVD}})\right| > D(\varepsilon)\right)$$

since
$$MASE(h_{\text{MASE}}) = \min_h MASE(h) \leq MASE(h_{\text{RCRVD}})$$

and
$$RCRVD(h_{\text{RCRVD}}) = \min_h RCRVD(h) \leq RCRVD(h_{\text{MASE}}).$$

CROSS-VALIDATION IN REGRESSION 19

Therefore, (A.13) can be written as

$$P(|h_{\text{RCRVD}} - h_{\text{MASE}}| > n^{-1/5}\varepsilon)$$
$$\leq P\left(n^{4/5}\sup_h \left|RCRVD(h_{\text{MASE}}) - \frac{E\psi'(\cdot)}{2}MASE(h_{\text{MASE}})\right| > \frac{D(\varepsilon)}{2}\right)$$
$$+ P\left(n^{4/5}\sup_h \left|\frac{E\psi'(\cdot)}{2}MASE(h_{\text{RCRVD}}) - RCRVD(h_{\text{RCRVD}})\right| > \frac{D(\varepsilon)}{2}\right)$$
$$\to 0,$$

by Proposition 1. Similarly, $P(|h_{\text{RCRVD}} - h_{\text{ASE}}| > n^{-1/5}\varepsilon) \to 0$ and $P(|h_{\text{MASE}} - h_{\text{ASE}}| > n^{-1/5}\varepsilon) \to 0$.  $\square$

**Acknowledgments.** The author wishes to thank Professor You Gan Wang and Professor Dawn J. Dekle for carefully reading through the paper. The author also wants to thank the editorial team for the many helpful comments and suggestions that have led to a substantially improved paper.

## REFERENCES

[1] AKAIKE, H. (1974). A new look at the statistical model identification. *IEEE Trans. Automatic Control* **19** 716–723. MR0423716
[2] BOENTE, G., FRAIMAN, R. and MELOCHE, J. (1997). Robust plug-in bandwidth estimators in nonparametric regression. *J. Statist. Plann. Inference* **57** 109–142. MR1440232
[3] CLEVELAND, W. S. (1979). Robust locally weighted regression and smoothing scatterplots. *J. Amer. Statist. Assoc.* **74** 829–836. MR0556476
[4] COX, D. D. (1983). Asymptotics for $M$-type smoothing splines. *Ann. Statist.* **11** 530–551. MR0696065
[5] GASSER, T. and MÜLLER, H.-G. (1979). Kernel estimation of regression functions. *Smoothing Techniques for Curve Estimation. Lecture Notes in Math.* **757** 23–68. Springer, Berlin. MR0564251
[6] HALL, P., SHEATHER, S., JONES, M. and MARRON, J. (1991). On optimal data-based bandwidth selection in kernel density estimation. *Biometrika* **78** 263–269. MR1131158
[7] HÄRDLE, W. (1984). How to determine the bandwidth of some nonlinear smoothers in practice. *Robust and Nonlinear Time Series Analysis. Lecture Notes in Statist.* **26** 163–184. Springer, New York. MR0786307
[8] HÄRDLE, W. (1987). Algorithm AS 222: Resistant smoothing using the fast Fourier transform. *Appl. Statist.* **36** 104–111.
[9] HÄRDLE, W. (1990). *Applied Nonparametric Regression*. Cambridge Univ. Press. MR1161622
[10] HÄRDLE, W. and GASSER, T. (1984). Robust non-parametric function fitting. *J. Roy. Statist. Soc. Ser. B* **46** 42–51. MR0745214
[11] HÄRDLE, W., HALL, P. and MARRON, J. (1988). How far are automatically chosen regression smoothing parameters from their optimum? (with discussion). *J. Amer. Statist. Assoc.* **83** 86–101. MR0941001

[12] HÄRDLE, W., HALL, P. and MARRON, J. (1992). Regression smoothing parameters that are not far from their optimum. *J. Amer. Statist. Assoc.* **87** 227–233. MR1158641
[13] HÄRDLE, W. and LUCKHAUS, S. (1984). Uniform consistency of a class of regression function estimators. *Ann. Statist.* **12** 612–623. MR0740915
[14] HÄRDLE, W. and MARRON, J. (1983). Optimal bandwidth selection in nonparametric function estimation. Institute of Statistics Mimeo Series no. 1530. Univ. North Carolina, Chapel Hill.
[15] HERRMANN, E. (2000). Variance estimation and bandwidth selection for kernel regression. In *Smoothing and Regression*: *Approaches, Computation and Application* (M. G. Schimek, ed.) 71–107. Wiley, New York.
[16] HUBER, P. (1964). Robust estimation of a location parameter. *Ann. Math. Statist.* **35** 73–101. MR0161415
[17] LEUNG, D., MARRIOTT, F. and WU, E. (1993). Bandwidth selection in robust smoothing. *J. Nonparametr. Statist.* **2** 333–339. MR1256384
[18] LIETZ, A. (2000). Analysis of water-quality trends at two discharge stations— one within Big Cypress National Preserve and one near Biscayne Bay—Southern Florida, 1966–1994. U.S. Geological Survey. Available at sofia.usgs.gov/publications/wri/00-4099.
[19] MARRON, J. and HÄRDLE, W. (1983). Random approximations to an error criterion of nonparametric statistics. Institute of Statistics Mimeo Series no. 1538. Univ. North Carolina, Chapel Hill.
[20] NADARAYA, E. (1964). On estimating regression. *Theory Probab. Appl.* **9** 141–142.
[21] PARK, B. and MARRON, J. (1990). Comparison of data-driven bandwidth selectors. *J. Amer. Statist. Assoc.* **85** 66–72.
[22] PRIESTLEY, M. and CHAO, M. (1972). Non-parametric function fitting. *J. Roy. Statist. Soc. Ser. B* **34** 385–392. MR0331616
[23] RICE, J. (1984). Bandwidth choice for nonparametric regression. *Ann. Statist.* **12** 1215–1230. MR0760684
[24] RONCHETTI, E., FIELD, C. and BLANCHARD, W. (1997). Robust linear model selection by cross-validation. *J. Amer. Statist. Assoc.* **92** 1017–1023. MR1482132
[25] SHIBATA, R. (1981). An optimal selection of regression variables. *Biometrika* **68** 45–54. MR0614940
[26] WANG, F. and SCOTT, D. (1994). The $L_1$ method for robust nonparametric regression. *J. Amer. Statist. Assoc.* **89** 65–76. MR1266287
[27] WATSON, G. (1964). Smooth regression analysis. *Sankhyā Ser. A* **26** 359–372. MR0185765
[28] WHITTLE, P. (1960). Bounds for the moments of linear and quadratic forms in independent variables. *Theory Probab. Appl.* **5** 302–305. MR0133849

SCHOOL OF ECONOMICS AND SOCIAL SCIENCES
SINGAPORE MANAGEMENT UNIVERSITY
90 STAMFORD ROAD
SINGAPORE 178903
REPUBLIC OF SINGAPORE
E-MAIL: leung@smu.edu.sg